\documentclass[12pt]{amsart}
\usepackage{ifthen,verbatim}
\usepackage{graphicx}
\usepackage{mathrsfs}
\usepackage{amsmath,amssymb,amsthm}
\usepackage{tikz-cd}
\usepackage{color}
\usepackage{caption}
\usepackage{subcaption}
\numberwithin{equation}{section}
\nonstopmode
\usepackage{amsfonts}
\usepackage{amssymb}
\usepackage[english]{babel}
\usepackage[utf8x]{inputenc}
\usepackage[autostyle]{csquotes}

\usepackage{float}
\usepackage{tabularx}
\newcolumntype{C}{>{\centering\arraybackslash}X}

\usepackage{tikz-cd}
\usepackage{color}
\usepackage{caption}
\newtheorem{remark}{Remark}
\newtheorem{theorem}{Theorem}
\newtheorem{lemma}{Lemma}

\newenvironment{myproof}[2] {\paragraph{\it Proof of {#1} {#2}.}}{\hfill$\square$}
\usepackage{url}
\usepackage{graphicx}
\usepackage{multirow}
\usepackage{tabularx}
\usepackage{longtable}
\usepackage{enumerate}

\usepackage{tikz-cd}

\usepackage{bm}
\usepackage{amsmath,mathtools}
\usepackage{amsfonts}

\newcommand{\RomanNumeralCaps}[1]
    {\MakeUppercase{\romannumeral #1}}

\usepackage[skip=1pt,font=scriptsize]{caption}

\newtheorem{innercustomthm}{Theorem}
\newenvironment{customthm}[1]
  {\renewcommand\theinnercustomthm{#1}\innercustomthm}
  {\endinnercustomthm}

\usepackage{dirtytalk}

\usepackage{mathtools}

\begin{document}
\title
{On Ramanujan's Modular Equations and Hecke Groups}

\makeatletter\def\thefootnote{\@arabic\c@footnote}\makeatother

\author[Md. S. Alam]{Md. Shafiul Alam}
\address{Department of Mathematics, University of Barishal, Barishal-8254, Bangladesh $-$ and $-$ Graduate School of Information Sciences,
Tohoku University, Aoba-ku, Sendai 980-8579, Japan}
\email{msalam@bu.ac.bd, shafiulmt@gmail.com}

\begin{abstract}
Inspired by the work of S. Ramanujan, many people have studied generalized modular equations and the numerous identities found by Ramanujan. These identities known as modular equations can be transformed into polynomial equations. There is no developed theory about how to find the degrees of these polynomial modular equations
explicitly. In this paper, we determine the degrees of the polynomial modular equations explicitly and study the relation between Hecke groups and modular equations in
Ramanujan’s theories of signatures 2, 3, and 4.


\end{abstract}
\thanks{The author was supported by Data Sciences Program \RomanNumeralCaps{2} of Graduate School of Information Sciences, Tohoku University, Japan.
}
\maketitle




\section{Introduction}\label{sec: Intro}
Let $\mathbb{D}$ denote the open unit disc $\{z\in \mathbb{C}:|z|<1\}$. For complex numbers $a,b,c$ with $c\neq 0,\,-1,\,-2,\,\dots$,  and nonnegative integer $n$, the Gaussian hypergeometric function, $_2F_1(a, b;c;z)$, is defined as
\begin{align*}
    _2F_1(a, b;c;z)=\sum_{n=0}^{\infty}\frac{(a)_n (b)_n}{(c)_n n!}z^n,\quad z\in \mathbb{D},
\end{align*}
where $(a)_n$ is the Pochhammer symbol or shifted factorial function given by
\begin{align*}
    (a)_n=
    \begin{cases}
      1, &  \text{if $n=0$}\\
      a(a+1)\cdots(a+n-1), & \text{if $n\geq1$.}
    \end{cases}
\end{align*}
By analytic continuation, $_2F_1(a, b;c;z)$ is extended to the slit plane $\mathbb{C}\setminus[1,\infty)$. For more details, see Chapter \RomanNumeralCaps{2} of \cite{bateman1953higher} and Chapter \RomanNumeralCaps{14} of \cite{whittaker2020course}.  

For $t\in(0,\,\frac{1}{2}],\,\alpha,\, \beta\in(0,1)$ and a given integer $p>1$, we say that $\beta$ has \emph{order} or \emph{degree} $p$ over $\alpha$ in the theory of signature $\frac{1}{t}$ if 
\begin{equation}\label{mod eq}
       \frac{_2F_1(t,1-t;1;1-\beta)}{_2F_1(t,1-t;1;\beta)}=p\,\frac{_2F_1(t,1-t;1;1-\alpha)}{_2F_1(t,1-t;1;\alpha)}.
    \end{equation}
Equation (\ref{mod eq}) is known as the generalized modular equation. In this article, we will use the terminology \emph{order} to avoid the confusion between the degree of the polynomial $P(\alpha,\beta)$ (see Theorem \ref{th:alam-sugawa21}) and the degree of the modulus $\beta$ over the modulus $\alpha$. The multiplier $m$ is given by \[m=\frac{{_2F_1(t,1-t;1;\alpha)}}{{_2F_1(t,1-t;1;\beta)}}.\] 
A modular equation of order $p$ in the theory of signature $\frac{1}{t}$ is an explicit relation between $\alpha$ and $\beta$ induced by (\ref{mod eq}) (see \cite{berndt1995ramanujan}). The great Indian mathematician S. Ramanujan extensively studied the generalized modular equation (\ref{mod eq}) and gave many identities involving $\alpha$ and $\beta$ for some rational values of $t$. Without original proofs, these identities were listed in Ramanujan's unpublished notebooks (see, e.g., \cite{berndt2012ramanujan}). There were no developed theories related to Ramanujan's modular equations before the 1980s. Some mathematicians, for example, B.~C. Berndt, S. Bhargava, J.~M. Borwein, P.~B. Borwein, F.~G. Garvan developed and organized the theories and tried to give the proofs of many identities recorded by S. Ramanujan (see \cite{berndt2012ramanujan, berndt1995ramanujan,borwein1987pi}). Also, G.~D. Anderson, M.~K. Vamanamurthy, M. Vuorinen and others have investigated the theory of Ramanujan's modular equations from different perspectives (see, e.g., \cite{anderson2000generalized,anderson1997conformal}).  

In this paper, we will consider the modular equations in the theories of signatures $2,\,3$, and $4$. There are different forms of modular equations for the same order of $\beta$ over $\alpha$ in the theory of signature $\frac{1}{t}$. For example, 
\begin{equation}\label{eq:me23}
(\alpha\beta)^{1/3}+\big\{(1-\alpha)(1-\beta)\big\}^{1/3}=1,
\end{equation}
\begin{equation}
     \Big\{\frac{(1-\beta)^2}{1-\alpha}\Big\}^\frac{1}{3}-\Big(\frac{\beta^2}{\alpha}\Big)^\frac{1}{3}=m
\end{equation}
and 
\begin{equation}
    \Big(\frac{\alpha^2}{\beta}\Big)^\frac{1}{3}+\Big\{\frac{(1-\alpha)^2}{1-\beta}\Big\}^\frac{1}{3}=\frac{4}{m^4}
\end{equation}
are the modular equations when the modulus $\beta$ has order $2$ over the modulus $\alpha$ in the theory of signature $3$ (see \cite[Theorem\, 7.1]{berndt1995ramanujan}). Note that \eqref{eq:me23} can be transformed to the following polynomial equation (see \cite{Alam-sugawa21})
\begin{equation*}
(2\alpha-1)^3\beta^3-3\alpha(4\alpha^2-13\alpha+10)\beta^2
+3\alpha(2\alpha^2-10\alpha+9)\beta-\alpha^3=0.
\end{equation*}

 

There is an intimate relation between the modular equations in Ramanujan's theories of signatures $\frac{1}{t}=2,\,3,\,4$ and the Hecke groups. The motivation of our present study comes from this relationship. The author and T. Sugawa \cite{Alam-sugawa21} offered a geometric approach to the proof of Ramanujan's identities for the solutions $(\alpha,\beta)$ to the generalized modular equation (\ref{mod eq}). They proved that the solution $(\alpha,\beta)$ satisfies a polynomial equation $P(\alpha,\beta)=0$. In this paper, we compute the degree in each of $\alpha$ and $\beta$ of the polynomial $P(\alpha,\beta)$ explicitly based on the relation between the Hecke groups and modular equations. We prove by geometric approach that if $(\alpha, \beta)$ is a solution to the generalized modular equation (\ref{mod eq}), then $(1-\beta,1-\alpha)$ is also a solution to (\ref{mod eq}) and $P(1-\beta,1-\alpha)=0$. Note that by the degree $\mu$ of the polynomial $P(\alpha, \beta)$, we will mean that  $P(\alpha,\beta)$ is a polynomial of degree $\mu$ in each of $\alpha$ and $\beta$.

For $t\in\Big\{\frac{1}{2},\,\frac{1}{3},\,\frac{1}{4}\Big\}$, let  
\begin{equation}\label{eq:lamda}
    \lambda_t=2\cos\frac{(1-2t)\pi}{2}
\end{equation}
and let $H(\lambda_t)$ denote the Hecke group generated by 
$$
A= \begin{pmatrix} 0 & -1\\ 1 & 0 \end{pmatrix}\quad\text{and}\quad B= \begin{pmatrix} 1 & \lambda_t\\ 0 & 1 \end{pmatrix}.
$$
If 
\begin{equation}\label{He}
   H_e(\lambda_t)=\Bigg\{\begin{pmatrix} a & b\lambda_t\\ c\lambda_t & d \end{pmatrix}: a,b,c,d \in\mathbb{Z} \text{ and } ad-bc\lambda_t^2=1\Bigg\}, 
\end{equation}
then $H_e(\lambda_t)$ is a subgroup of $H(\lambda_t)$ of index $2$ (see \cite{cangul1998normal}). Note that $H_e(\lambda_t)$ is called the even subgroup of $H(\lambda_t)$ for $\lambda_t=\sqrt{2}$ and $\sqrt{3}$. We will consider $H_e(\lambda_t)$ for $\lambda_t=\sqrt{2}, \sqrt{3}$ and $2$. Let $\mathbb{H}$ denote the upper half-plane $\{\tau\in \mathbb{C}: \text{Im } \tau>0\}$. Then the quotient Riemann surface $H_e(\lambda_t)\backslash\mathbb{H}$ is $\widehat{\mathbb{C}}\setminus\{0,1\}$ for $t\in\Big\{\frac{1}{3},\,\frac{1}{4}\Big\}$ and $\widehat{\mathbb{C}}\setminus\{0,1,\infty\}$ for $t=\frac{1}{2}$. The following theorem asserts that the solution $(\alpha,\beta)$ to the generalized modular equation (\ref{mod eq}) satisfies a polynomial equation in $\alpha$ and $\beta$.

\begin{customthm}{A}[$\text{\cite[Theorem 1.8]{Alam-sugawa21}}$]\label{th:alam-sugawa21}
For integers $p, n>1$ and $t\in(0,1/2]$, let \[H_e'(\lambda_t)=M_p^{-1}H_e(\lambda_t)M_p \quad \text{and}\quad H_{M_p}(\lambda_t)=H_e(\lambda_t)\cap H_e'(\lambda_t),\] 
where $M_p=\begin{pmatrix} p & 0\\ 0 & 1 \end{pmatrix}$. Then, the solution $(\alpha, \beta)$ to the generalized modular equation (\ref{mod eq}) in $H_e(\lambda_t)\backslash\mathbb{H}$ satisfies the polynomial equation $P(\alpha,\beta)=0$ for an irreducible polynomial $P(x,y)$ of degree $\mu$ in each of $x$ and $y$ if and only if $H_{M_p}(\lambda_t)$ is a subgroup of $H_e(\lambda_t)$ of index $\mu$.
\end{customthm}

H. H. Chan and W.-C. Liaw \cite{chan2000russell} studied modular equations in the theory of signature $3$ based on the modular equations studied by R. Russell \cite{Russell}.

\begin{customthm}{B}[$\text{\cite[Theorems 2.1, 3.1]{chan2000russell}}$]\label{th:russell}

If $p>2$ is a prime, $u=(\alpha\beta)^{l/8}$ and $v=\{(1-\alpha)(1-\beta)\}^{l/8}$, where $(p+1)/8=m/l$ in lowest terms, then $(u,v)$ satisfies a polynomial equation $Q(u,v)=0$, where $Q(x,y)$ is of degree $m$ in each of $x$ and $y$ in the theory of signature $2$. If $p>3$ is a prime, $u=(\alpha\beta)^{l/6}$ and $v=\{(1-\alpha)(1-\beta)\}^{l/6}$, where $(p+1)/3=m/l$ in lowest terms, then $(u,v)$ satisfies a polynomial equation $Q(u,v)=0$, where $Q(x,y)$ is of degree $m$ in each of $x$ and $y$ in the theory of signature $3$.

\end{customthm}

\begin{remark}
In the theory of signature $3$, the degree $\mu$ of the polynomial $P(\alpha, \beta)$ in Theorem \ref{th:alam-sugawa21} and the degree $m$ of the polynomial $Q(u,v)$ in Theorem \ref{th:russell} are related as follows:
\begin{enumerate}[(i)]
\item $\mu=3m$ when $p\equiv 2\,(\text{mod }3)$,
\item $\mu=m$ when $p\equiv 1\,(\text{mod }3)$.

\end{enumerate}

\end{remark}

The remainder of this article is organized as follows. In Section \ref{sec:main results}, we state our main results. Some basic facts related to modular groups and Hecke groups are discussed in Section \ref{sec:preli}. Finally, the proofs of the main results are given in Section \ref{sec:proofs}.

\section{Main Results}\label{sec:main results}
Let $\Psi(N)$ denote the Dedekind psi function given by 
\begin{equation}\label{eq:psi}
    \Psi(N)=N\prod_{\underset{q \text{ prime}}{q|N}}\Big(1+\frac{1}{q}\Big),\quad N\in\mathbb{N}
\end{equation}
(see \cite[p.\,123]{dickson}). Our first result is for determining the degree in each of $\alpha$ and $\beta$ of the polynomial $P(\alpha,\beta)$ in Theorem \ref{th:alam-sugawa21} explicitly in Ramanujan's theories of signatures $\frac{1}{t}=2,\,3$ and $4$. 

\begin{theorem}\label{theorem degree}
For an integer $p>1$, suppose $\beta$ has order $p$ over $\alpha$ in the theories of signatures $\frac{1}{t}=2$, $3$ and $4$. Let $\mu\big(p,\frac{1}{t}\big)$ be the degree in each of $\alpha$ and $\beta$ of the polynomial $P(\alpha,\beta)$, then
$$
\mu(p,2)=\mu(p,4)=\frac{1}{3}\Psi(2p)
$$
and 
$$\mu(p,3)=\frac{1}{4}\Psi(3p).$$

\end{theorem}

\begin{remark}
If $p$ is an odd prime, then $\mu(p,2)=\mu(p,4)=p+1$. If $p\neq3$ is a prime, then $\mu(p,3)=p+1$.
\end{remark}

We compute the degree $\mu\big(p,\frac{1}{t}\big)$ for some small values of $p$ and $t\in\Big\{\frac{1}{2},\,\frac{1}{3},\,\frac{1}{4}\Big\}$ in Table \ref{tab: index}. Even if one does not know the corresponding Hecke subgroups, he/she can compute the degree of modular equations in the theories of signatures 2, 3, and 4 using the formulas in Theorem \ref{theorem degree}.

The following result establishes some statements related to the Hecke subgroups and the modular equations in the theories of signatures 2, 3, and 4.

\begin{theorem}\label{theorem orbit}
For a given integer $p>1$, suppose that $\beta$ has order $p$ over $\alpha$ in the theories of signatures $\displaystyle\frac{1}{t}=2,\,3$ and $4$. 
If $M_p=\begin{pmatrix} p & 0\\ 0 & 1 \end{pmatrix}$, then 

\begin{enumerate}[(i)]
    \item there exists a Hecke subgroup, say $\Gamma_1$, of finite index in $H(\lambda_t)$,
    \item $\big(M_p^{-1}\Gamma_1M_p\big)\cap\Gamma_1$ has finite index in $H(\lambda_t)$,
    \item the degree of the branched covering \[\Gamma_2\backslash\mathbb{H}\rightarrow \Gamma_1\backslash\mathbb{H}\]
    is finite, where $\Gamma_2=\big(M_p^{-1}\Gamma_1M_p\big)\cap\Gamma_1$,
    \item there is a polynomial equation $P(\alpha,\beta)=0$, where the polynomial $P(\alpha,\beta)$ has degree $\mu=\big|\Gamma_1:\Gamma_2\big|$ in each of $\alpha$ and $\beta$.
\end{enumerate}
\end{theorem}

\begin{remark}
In fact, the statements in Theorem \ref{theorem orbit} are mutually equivalent.
\end{remark}

\begin{table}[h]

\begin{equation*}
\begin{tabular}{ |c|c|c|c| } 
\hline
 & & \\[-.5em] 
$p$ & $\mu(p,2)$ and $\mu(p,4)$ & $\mu(p,3)$\\& &\\[-.5em] 

\hline
2 & 2 & 3\\ 
3 & 4 & 3\\ 
4 & 4 & 6\\ 
5 & 6 & 6\\ 
6 & 8 & 9\\ 
7 & 8 & 8\\ 
8 & 8 & 12\\
9 & 12 & 9\\ 
10 & 12 & 18\\ 
11 & 12 & 12\\ 
12 & 16 & 18\\ 
13 & 14 & 14\\ 
14 &  16 & 24\\ 
15 & 24 & 18\\ 
16 & 16 & 24\\
17 & 18 & 18\\ 
18 & 24 & 27\\
19 & 20 & 20\\ 
20 & 24 & 36\\
\hline

\end{tabular}
\end{equation*}
\caption{\label{tab: index} Values of $\mu\big(p,\frac{1}{t}\big)$ for some small values of $p$ and $t\in\Big\{\frac{1}{2},\,\frac{1}{3},\,\frac{1}{4}\Big\}$.}
\end{table}

We can express the generalized modular equation (\ref{mod eq}) as $f_t(\beta)=p f_t(\alpha)$, where $f_t$ is defined by
\begin{equation}\label{map:ft}
   f_t(z)=i\,\frac{_2F_1(t,1-t;1;1;1-z)}{_2F_1(t,1-t;1;1;z)}, 
\end{equation}
$t\in(0,1/2]$ and $p$ is an integer $>1$. Consider the canonical projection $\pi_t:\mathbb{H}\rightarrow H_e(\lambda_t)\backslash\mathbb{H}$, where $\pi_t$ is the inverse of $f_t$ (a detailed discussion will be given in Section \ref{sec:preli}). The moduli $\alpha, \beta\in\widehat{\mathbb{C}}\setminus\{0,1\}$ satisfy (\ref{mod eq}) if and only if $\alpha=\pi_t(\tau)$ and $\beta(\tau)=\pi_t(p\tau)$ for $\tau\in\mathbb{H}$ (see \cite{Alam-sugawa21}) and we have the following theorem.

\begin{theorem}\label{th:alpha,beta}
For the canonical projection $\pi_t:\mathbb{H}\rightarrow H_e(\lambda_t)\backslash\mathbb{H}$, let $\alpha=\pi_t(\tau)$ and $\beta=\pi_t(p\tau)$, where $\tau\in\mathbb{H}$ and $p$ is an integer $>1$. If the solution $(\alpha,\beta)$ to the generalized modular equation (\ref{mod eq}) satisfies the equation $P(x,y)=0$, then $(1-\beta,1-\alpha)$ is also a solution to (\ref{mod eq}) and satisfies the equation  $P(x,y)=0$, where $P(x,y)$ is the polynomial in Theorem \ref{th:alam-sugawa21}. 
\end{theorem}

\begin{remark}
    The equation $P(1-\beta,1-\alpha)=0$ is the reciprocal of the equation $P(\alpha, \beta)=0$ and this process is known as the method of reciprocation (see Theorem 6.3.2 in \cite{berndt2006}). We use Lemma 3.4 of \cite{Alam-sugawa21} in the proof of Theorem 3 and our proof is geometric. 

\end{remark}

\vspace{5mm}

\section{Preliminaries}\label{sec:preli}


The group $\operatorname{SL}(2,\mathbb{R})$ is defined by
$$ 
\operatorname{SL}(2,\mathbb{R})=\Bigg\{\displaystyle\begin{pmatrix} a & b\\ c & d \end{pmatrix}: a,b,c,d\in\mathbb{R},ad-bc=1\Bigg\}
$$
and is generated by
\[ \begin{pmatrix} 0 & -1\\ 1 & 0 \end{pmatrix}\quad\text{and}\quad \begin{pmatrix} 1 & 1\\ 0 & 1 \end{pmatrix}.\]
Let $I_2$ denote the $2\times2$ identity matrix, then the group $\operatorname{PSL}(2,\mathbb{R})=\operatorname{SL}(2, \mathbb{R})/\{\pm I_2\}$. For $a, b, c, d\in \mathbb{R}$ and $ad-bc=1$, the group $\operatorname{PSL}(2,\mathbb{R})$ acts on the upper half-plane $\mathbb{H}$ as follows
   \[\tau\mapsto \frac{a\tau+b}{c\tau+d}\]
and $\operatorname{PSL}(2,\mathbb{R})$ is the group of automorphisms of the upper half-plane $\mathbb{H}$. All transformations of $\operatorname{PSL}(2,\mathbb{R})$ are conformal.  Assume that $\Gamma$ is a Fuchsian group of the first kind which leaves the upper half-plane $\mathbb{H}$ or the unit disc $\mathbb{D}$ invariant. Then, $\Gamma$ is a discrete subgroup of the group of orientation-preserving isometries of $\mathbb{H}$, i.e., $\Gamma$ is a discrete subgroup of $\operatorname{PSL}(2, \mathbb{R})$ (see \cite{katok1992fuchsian}). 

If $m$ is a positive integer, then the congruence subgroup $\Gamma_0(m)$ is defined as 
\begin{equation*}
    \Gamma_0(m)=\Bigg\{\begin{pmatrix} a & b\\ c & d \end{pmatrix}\in \operatorname{PSL}(2,\mathbb{Z}): c\equiv0\hspace{.5mm}(\text{mod}\hspace{1mm}m)\Bigg\}.
\end{equation*}

We now construct a connection between the Schwarz triangle function and the Gaussian hypergeometric function, $_2F_1(a, b;c;z)$. Since \[w_1={_2F_1}(a, b;c;z)\] and \[w_2={_2F_1}(a, b;a+b+1-c;1-z)\] are two linearly independent solutions of the following hypergeometric differential equation
\begin{equation}\label{hyp}
    z(1-z)\frac{d^2w}{dz^2}+\{c-(a+b+1)z\}\frac{dw}{dz}-abw=0,
\end{equation}
it is a well-known fact that the Schwarz triangle function defined by
\begin{equation*}
    S(z)=i\,\frac{{_2F_1}(a, b;a+b+1-c;1-z)}{{_2F_1}(a, b;c;z)}
\end{equation*}
maps the upper half-plane $\mathbb{H}$ conformally onto a curvilinear triangle $\Delta_t$, which has interior angles $(1-c)\pi$, $(c-a-b)\pi$ and $(b-a)\pi$ at the vertices $S(0)$, $S(1)$ and $S(\infty)$, respectively. For details, we recommend the readers to go through Chapter \RomanNumeralCaps{5}, Section 7 of \cite{nehari1952conformal}. 
For $t\in(0,\frac{1}{2}]$, let $a=t,\, b=1-a=1-t$ and $c=1$, then $S(z)$ can be expressed as
\begin{equation}\label{map:ft}
   S(z)=f_t(z)=i\,\frac{_2F_1(t,1-t;1;1;1-z)}{_2F_1(t,1-t;1;1;z)}. 
\end{equation}

If $\theta_1=\frac{\pi}{m_1},\  \theta_2=\frac{\pi}{m_2}$ and $\theta_3=\frac{\pi}{m_3}$, then a curvilinear triangle with angles $\theta_j$ (for $j=1,\ 2,\ 3$) can be continued as a single-valued function across the sides of the triangle $\Delta_t$ if and only if $m_j$ is an integer greater than $1$ including $\infty$ (see  \cite[p.~416]{gerretsen1969geometric}). Therefore,
\begin{equation}\label{1/m}
    \frac{1}{m_1}+\frac{1}{m_2}+\frac{1}{m_3}<1.
\end{equation}
The following lemma is related to the above facts.
\begin{lemma}[$\text{\cite[Lemma 4.1]{anderson2010twice}}$]\label{lem: H to delta} Let the map $f_t$ be defined by (\ref{map:ft}) for $t\in(0,\frac{1}{2}]$. Then, the upper half-plane $\mathbb{H}$ is mapped by $f_t$ onto the hyperbolic triangle $\Delta_t$ given by 

\[\Delta_t=\Bigg\{\tau\in\mathbb{H}:0< Re\hspace{.8mm}\tau<\cos\frac{\theta}{2}, \Big|2\tau\cos\frac{\theta}{2}-1\Big|>1\Bigg\},\] 
where $\theta=(1-2t)\pi$. The interior angles of $\Delta_t$ are $0,0,$ and $\theta=(1-2t)\pi$ at the vertices $f_t(0)=i\infty$, $f_t(1)=0$ and $f_t(\infty)=e^{i\frac{\theta}{2}}$, respectively.  
\end{lemma}

\begin{figure}[h]
         \centering
         \includegraphics[width=0.8\textwidth]{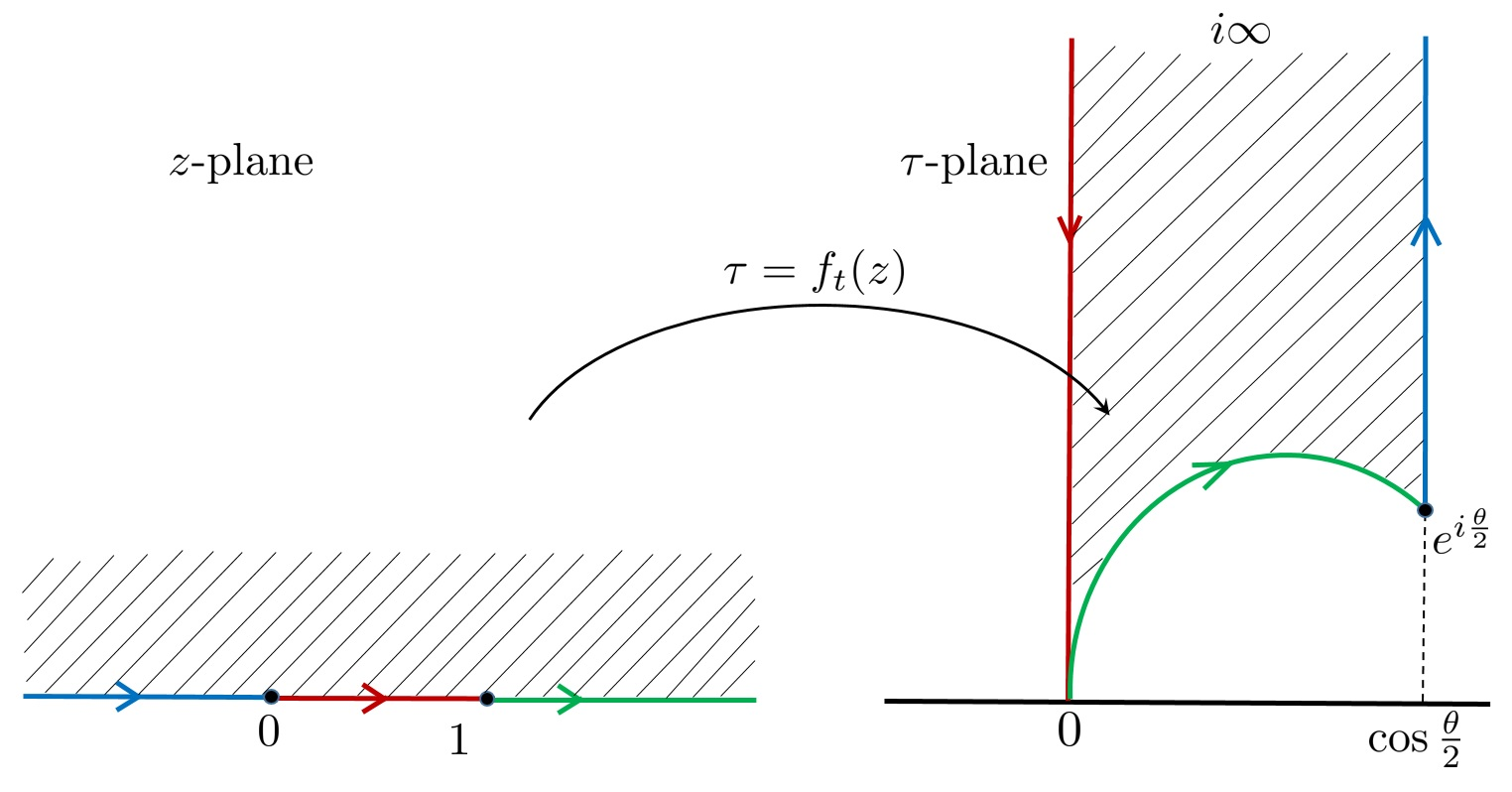}
         \caption{Mapping of the upper half-plane $\mathbb{H}$ onto $\Delta_t$ by $f_t$.}
         \label{fig:32}
\end{figure}



By Lemma \ref{lem: H to delta}, the condition (\ref{1/m}) becomes $\frac{1}{m_3}<1$,
i.e., it depends only on the third fixed point $f_t(\infty)=e^{i\frac{\theta}{2}}$ and $m_3=\frac{1}{1-2t}$ is an integer greater than $1$ including $\infty$ only for $t\in\Big\{\frac{1}{2}, \frac{1}{3},\frac{1}{4}\Big\}$. If $\pi_t: \Delta_t\rightarrow \mathbb{H}$ is the inverse map of $f_t$, then we can extend $\pi_t(\tau)$ analytically to a single-valued function on $\mathbb{H}$ with the real axis as its natural boundary by applying the Schwarz reflection principle repeatedly. The covering group of $\pi_t$ is the Hecke subgroup $H_e(\lambda_t)$. For more details, see Section 2 of \cite{Alam-sugawa21}, where $H_e(\lambda_t)$ is denoted by $G_q$.

The subgroup $H_e(\lambda_t)$ has two cusps and one elliptic point for $t\in\Big\{\frac{1}{3},\,\frac{1}{4}\Big\}$ and has three cusps for $t=\frac{1}{2}$. Thus, the quotient Riemann surface $H_e(\lambda_t)\backslash\mathbb{H}$ is the two punctured Riemann sphere $\widehat{\mathbb{C}}\setminus\{0,1\}$ for $t\in\Big\{\frac{1}{3},\,\frac{1}{4}\Big\}$ and the thrice punctured Riemann sphere $\widehat{\mathbb{C}}\setminus\{0,1,\infty\}$ for $t=\frac{1}{2}$. The set of cusps of the Hecke group $H(\lambda_t)$ is $\mathbb{Q}[\lambda_t]\cup\{\infty\}$. To compactify the quotient Riemann surface $H_e(\lambda_t)\backslash\mathbb{H}$, let $\mathbb{H}^*=\mathbb{H}\cup\mathbb{Q}[\lambda_t]\cup\{\infty\}$. Then, $H_e(\lambda_t)\backslash\mathbb{H}^*$ is a compact Riemann surface. For all $\begin{pmatrix} a & b\\ c & d \end{pmatrix}\in H(\lambda_t)$ and $\tau\in\mathbb{H}$, the meromorphic function $g:\mathbb{H}\rightarrow H(\lambda_t)\backslash\mathbb{H}^*$ is called an automorphic function if $g\bigg(\displaystyle\frac{a\tau+b}{c\tau+d}\bigg)=g(\tau)$ (see \cite{Bump}).

\section{Proofs of Main Results}\label{sec:proofs}

Let $\Gamma=\operatorname{PSL}(2,\mathbb{Z})$. For an integer $p>1$, let $M_p=\begin{pmatrix} p & 0\\ 0 & 1 \end{pmatrix}$, then the transformation group of order $p$ (see Chapter \RomanNumeralCaps{6} of \cite{schoeneberg1974elliptic}), $\Gamma_{M_p}$,  is given by
\[\Gamma_{M_p}:=\Gamma\cap \Big(M_p^{-1}\Gamma M_p\Big),\]
which can be written as the group of M\"obius transformations
\[\Gamma_{M_p}:=\Big\{\gamma\in\Gamma:\hspace{1mm} M_p\gamma M_p^{-1}\in\Gamma\Big\}.\]
If $\gamma=\begin{pmatrix} a & b\\ c & d \end{pmatrix}\in\Gamma$, then $M_p\gamma M_p^{-1}=\begin{pmatrix} a & p\hspace{.3mm}b\\ \frac{c}{p} & d \end{pmatrix}$. Hence, $M_p\gamma M_p^{-1}\in\Gamma$ only when $c\equiv 0\,(\text{mod}\,p)$, and we have $\Gamma_{M_p}=\Gamma_0(p)$.
The following lemma is a well-known result, e.g., see Proposition 1.43 in \cite{shimura1971introduction} or  \cite[ p.\,79]{schoeneberg1974elliptic}.

\begin{lemma}\label{lem index}
For any positive integer $N$, $\big|\Gamma:\Gamma_0(N)\big|=\Psi(N)$.
\end{lemma}

\begin{myproof}{Theorem}{\ref{theorem degree}}
For $t\in\Big\{\frac{1}{2},\,\frac{1}{3},\,\frac{1}{4}\Big\}$ and $\lambda_t=2\cos\frac{(1-2t)\pi}{2}$, let
\[H_{M_p}(\lambda_t)=\Big\{\gamma\in H_e(\lambda_t):M_p\gamma M_p^{-1}\in H_e(\lambda_t)\Big\},\]
where $p$ is an integer $>1$ and $M_p=\begin{pmatrix} p & 0\\ 0 & 1 \end{pmatrix}$. If $\gamma=\begin{pmatrix} a & b\lambda_t\\ c\lambda_t & d \end{pmatrix}\in H_e(\lambda_t)$, then 
$$
M_p\gamma M_p^{-1}=\begin{pmatrix} a & p\hspace{.5mm}b\lambda_t\\ \frac{c}{p}\hspace{.5mm}\lambda_t & d \end{pmatrix}.
$$ 
Therefore, $M_p\gamma M_p^{-1}\in H_e(\lambda_t)$ only when $c\equiv 0\,(\text{mod }p)$ and we have
\begin{equation}\label{HM_p}
   H_{M_p}(\lambda_t)=\Bigg\{\begin{pmatrix} a & b\lambda_t\\ c\lambda_t & d \end{pmatrix}\in H_e(\lambda_t): c\equiv 0\,(\text{mod }p)\Bigg\}. 
\end{equation}
Consequently,
\begin{equation}\label{itrsec of sg}
    H_e(\lambda_t)\cap \Big(M_p^{-1}H_e(\lambda_t)M_p\Big)=H_{M_p}(\lambda_t)
\end{equation}
and \[H_{M_p}(\lambda_t)<H_e(\lambda_t)<H(\lambda_t).\]
Let $\pi_t$ and $\pi_t'$ denote the canonical projections $\mathbb{H}\rightarrow H_e(\lambda_t)\backslash\mathbb{H}$ and $\mathbb{H}\rightarrow H_{M_p}(\lambda_t)\backslash\mathbb{H}$, respectively. From the subgroup relation $H_{M_p}(\lambda_t)<H_e(\lambda_t)$, we have the branched covering map $\varphi:H_{M_p}(\lambda_t)\backslash\mathbb{H}\rightarrow H_e(\lambda_t)\backslash\mathbb{H}$ and the following commutative diagram:
\begin{center}
   \begin{tikzcd}[column sep=40pt,row sep=40pt]
&\mathbb{H} \arrow[d, "\pi_t'", labels=left] \arrow[dr, "\pi_t"] & \\
&H_{M_p}(\lambda_t)\backslash\mathbb{H} \arrow[r, "\varphi", labels=below]  &H_e(\lambda_t)\backslash\mathbb{H}.
\end{tikzcd}
\end{center}
The degree of the branched covering $H_{M_p}(\lambda_t)\backslash\mathbb{H}\rightarrow H_e(\lambda_t)\backslash\mathbb{H}$ is $\big|H_e(\lambda_t):H_{M_p}(\lambda_t)\big|$, which is the degree $\mu\big(p,\frac{1}{t}\big)$ in each of $\alpha$ and $\beta$ of the polynomial $P(\alpha,\beta)$ by Theorem \ref{th:alam-sugawa21}.

Also, for $\Gamma=\operatorname{PSL}(2,\mathbb{Z})$, we have
\[\Gamma\cap \Big(M_p^{-1}\Gamma M_p\Big)=\Gamma_0(p),\]
 and 
\[\Gamma_0(\lambda_t^2p) < \Gamma_0(\lambda_t^2) < \Gamma.\]

Let us consider the mapping \[\Theta:H_e(\lambda_t) \rightarrow \Gamma_0(\lambda_t^2)\] defined by
\[\Theta(A)=M_{\lambda_t}^{-1}AM_{\lambda_t},\]
where $A=\begin{pmatrix} a & b\lambda_t\\ c\lambda_t & d \end{pmatrix}\in H_e(\lambda_t)$ and $M_{\lambda_t}=\begin{pmatrix} \lambda_t & 0 \\ 0 & 1 \end{pmatrix}$. Then, we have
\begin{align*}
   \Theta(H_e(\lambda_t))= \Gamma_0(\lambda_t^2)\quad \text{and}\quad
    \Theta(H_{M_p}(\lambda_t))= \Gamma_0(\lambda_t^2p).
\end{align*}
Therefore, $H_e(\lambda_t)\cong\Gamma_0(\lambda_t^2)$, $H_{M_p}(\lambda_t)\cong\Gamma_0(\lambda_t^2p)$, and we have
\begin{align*}
     \big|H_e(\lambda_t):H_{M_p}(\lambda_t)\big|&=\big|\Gamma_0(\lambda_t^2):\Gamma_0(\lambda_t^2p)\big|=\frac{\big|\Gamma:\Gamma_0(\lambda_t^2p)\big|}{\big|\Gamma:\Gamma_0(\lambda_t^2)\big|}.
\end{align*}
By Lemma \ref{lem index}, $\big|\Gamma:\Gamma_0(\lambda_t^2p)\big|=\Psi(\lambda_t^2p)$ and $\big|\Gamma:\Gamma_0(\lambda_t^2)\big|=\Psi(\lambda_t^2)$. Hence $$\mu\Big(p,\frac{1}{t}\Big)=\big|H_e(\lambda_t):H_{M_p}(\lambda_t)\big|=\frac{\Psi(\lambda_t^2p)}{\Psi(\lambda_t^2)},$$ 
which implies $\mu(p,2)=\frac{1}{6}\Psi(4p)$, $\mu(p,3)=\frac{1}{4}\Psi(3p)$ and $\mu(p,4)=\frac{1}{3}\Psi(2p)$. By (\ref{eq:psi}), it is easy to show that $\Psi(4p)=2\Psi(2p)$. Thus, $\mu(p,2)=\mu(p,4)$ as required.
\end{myproof} 

\vspace{5mm}


Let $X_1=H_e(\lambda_t)\backslash\mathbb{H}$ and $X_2=H_{M_p}(\lambda_t)\backslash\mathbb{H}$. For the canonical projections $\pi_t:\mathbb{H}\rightarrow X_1$ and $\pi_t':\mathbb{H}\rightarrow X_2$, consider the mappings $\varphi:X_2\rightarrow X_1$ and $\psi:X_2\rightarrow X_1$ such that $\pi_t=\varphi\circ\pi_t'$ and $\pi_t\circ M_p=\psi\circ\pi_t'$, i.e., the following diagrams commute:
\begin{equation*}
\begin{tikzcd}
&\mathbb{H} \arrow[d, "\pi_t'", labels=left] \arrow[dr, "\pi_t"] & \\
& X_2 \arrow[r, "\varphi", labels=below]  &X_1
\end{tikzcd}\hspace{40pt}%
\begin{tikzcd}
&\mathbb{H} \arrow[r, "M_p"] \arrow[d, "\pi_t'",labels=left]
& \mathbb{H} \arrow[d, "\pi_t"] \\
&X_2 \arrow[r, "\psi",labels=below]
& X_1.
\end{tikzcd}
\end{equation*}
Thus, for $z\in X_2$, the solution $(\alpha, \beta)$ to the generalized modular equation (\ref{mod eq}) is parametrized by $\alpha=\varphi(z)$ and $\beta=\psi(z)$. Before giving the proofs of Theorems \ref{theorem orbit} and \ref{th:alpha,beta}, we recall the following two lemmas from \cite{Alam-sugawa21}.

\begin{lemma}[$\text{\cite[Lemma 2.3]{Alam-sugawa21}}$]\label{lem index alam-sugawa}
For an integer $p>1$, let $t\in(0,\frac{1}{2}]$ and $M_p=\begin{pmatrix} p & 0\\ 0 & 1 \end{pmatrix}$. If $H_e'(\lambda_t)=M_p^{-1}H_e(\lambda_t)M_p$ and $H_{M_p}(\lambda_t)=H_e(\lambda_t)\cap H_e'(\lambda_t)$, then 
$$
\big|H_e(\lambda_t):H_{M_p}(\lambda_t)\big|=\big|H_e'(\lambda_t):H_{M_p}(\lambda_t)\big|.
$$
\end{lemma}

\begin{lemma}[$\text{\cite[Lemma 3.4]{Alam-sugawa21}}$]\label{lem:Fricke}
Let $W_p= \begin{pmatrix} 0 & -1\\ p & 0 \end{pmatrix}$ and $\tau\in\mathbb{H}$, then $W_p\tau=\displaystyle-\frac{1}{p\tau}$ induces an automorphism $\omega$ on $X_2$ such that $\pi_t'\circ W_p=\omega\circ\pi_t'$ and $\psi=1-\varphi\circ\omega$.
\end{lemma}

Now, we are ready to prove Theorems \ref{theorem orbit} and \ref{th:alpha,beta}.

\hspace{.5cm}
\begin{myproof}{Theorem}{\ref{theorem orbit}} 
First, recall that the covering group of the map $\pi_t$ is the Hecke subgroup $H_e(\lambda_t)$ and it is well-known that the index of $H_e(\lambda_t)$ in the Hecke group $H(\lambda_t)$ is 2 (see \cite[p. 61]{cangul1998normal}). Thus, $\Gamma_1=H_e(\lambda_t)$ and $(\romannumeral 1)$ follows easily from this fact.

From (\ref{itrsec of sg}), we have $\Gamma_2=\Gamma_1\cap \big(M_p^{-1}\Gamma_1M_p\big)=H_{M_p}(\lambda_t)$. By virtue of the proof of Theorem \ref{theorem degree}, we have $\Gamma_2\cong\Gamma_0(\lambda_t^2p)$, and hence, $\Gamma_2$ is isomorphic to $\Gamma_0(4p)$, $\Gamma_0(3p)$ and $\Gamma_0(2p)$ for $t=\frac{1}{2}, \frac{1}{3}$ and $\frac{1}{4}$, respectively. Each of $\Gamma_0(4p)$, $\Gamma_0(3p)$ and $\Gamma_0(2p)$ has finite index in $\Gamma=\operatorname{PSL}(2,\mathbb{Z})$. Therefore, $\Gamma_2$ has finite index in $H(\lambda_t)$, which implies $(\romannumeral 2)$.

If $X_1=\Gamma_1\backslash\mathbb{H}$ and $X_2=\Gamma_2\backslash\mathbb{H}$, then the degree of the branched covering $\varphi:X_2\rightarrow X_1$ is equal to the index of $\Gamma_2$ in $\Gamma_1$. Since each of $\Gamma_1$ and $\Gamma_2$ has finite index in $H(\lambda_t)$, the index of $\Gamma_2$ in $\Gamma_1$ is finite. Therefore, $(\romannumeral 3)$ holds.

It is not difficult to prove that $\alpha(\tau)$ and $\beta(\tau)=\alpha(p\tau)$ are automorphic functions on $\Gamma_1$ and $\Gamma_1':=M_p^{-1}\Gamma_1M_p$, respectively. Recall that the quotient Riemann surface $X_1=\Gamma_1\backslash\mathbb{H}$ is $\widehat{\mathbb{C}}\setminus\{0,1\}$ for $t\in\Big\{\frac{1}{3},\,\frac{1}{4}\Big\}$ and $\widehat{\mathbb{C}}\setminus\{0,1,\infty\}$ for $t=\frac{1}{2}$. If $\widehat{X}_1$ is the compactification of $X_1$, then $\widehat{X}_1$ is the Riemann sphere $\widehat{\mathbb{C}}$. Thus, the field of automorphic functions for $\Gamma_1$ is $\mathbb{C}(\alpha(\tau))$. If $X_1'=\Gamma_1'\backslash\mathbb{H}$ and $\widehat{X}_1'$ is the compactification of $X_1'$, then $\widehat{X}_1'=\widehat{\mathbb{C}}$. The field of automorphic functions for $\Gamma_1'$ is $\mathbb{C}(\beta(\tau))$.  Since $\Gamma_2<\Gamma_1$ and $\Gamma_2<\Gamma_1'$, both $\mathbb{C}(\alpha(\tau))$ and $\mathbb{C}(\beta(\tau))$ are subfields of the field of automorphic functions for $\Gamma_2=\Gamma_1\cap\Gamma_1'$, i.e., $\mathbb{C}(\alpha(\tau), \beta(\tau))$. If $\mu=\big|\Gamma_1:\Gamma_2\big|$, then $\varphi: X_2\rightarrow X_1$ is a $\mu$-sheeted branched covering map. For any function $g\in\mathbb{C}(\alpha(\tau))$, we have a function $f\in\mathbb{C}(\alpha(\tau), \beta(\tau))$ by virtue of the pullback $\varphi^*(g)=g\circ \varphi=f$, where $\varphi^*: \mathbb{C}(\alpha(\tau))\rightarrow \mathbb{C}(\alpha(\tau), \beta(\tau))$ is an algebraic field extension of degree $\mu$ (see \cite[Theorem 8.3]{forster2012lectures}). Similarly, if $\psi$ is the branched covering map $X_2\rightarrow X_1'$, then $\psi$ is also a $\mu$-sheeted covering map, since $|\Gamma_1':\Gamma_2|=\mu$ by Lemma \ref{lem index alam-sugawa}. Hence, $\psi^*: \mathbb{C}(\beta(\tau))\rightarrow \mathbb{C}(\alpha(\tau), \beta(\tau))$ is an algebraic field extension of degree $\mu$. Consequently, there is a polynomial $P(\alpha(\tau),\beta(\tau))$ which has degree $\mu$ in each of $\alpha(\tau)$ and $\beta(\tau)$. The polynomial $P(\alpha(\tau),\beta(\tau))$ is determined up to a scalar factor so that $P(\alpha(\tau),\beta(\tau))=0$, which implies $(\romannumeral 4)$ and completes the proof.
\end{myproof}

\hspace{.5cm}

Recall that the Hecke subgroup $H_{M_p}(\lambda_t)$ is given by 
\[H_{M_p}(\lambda_t)=\Bigg\{\begin{pmatrix} a & b\lambda_t\\ c\lambda_t & d \end{pmatrix}\in H_e(\lambda_t): c\equiv 0\,(\text{mod }p)\Bigg\}.\]

Let $W_p= \begin{pmatrix} 0 & -1\\ p & 0 \end{pmatrix}$, then
\begin{equation}\label{w_p thm3}
    W_p^{-1}\begin{pmatrix} a & b\lambda_t\\ c\lambda_t & d \end{pmatrix}W_p=\begin{pmatrix} d & -\frac{c}{p}\lambda_t\\ -pb\lambda_t & a\end{pmatrix},
\end{equation}
where $\begin{pmatrix} a & b\lambda_t\\ c\lambda_t & d \end{pmatrix}\in H_{M_p}(\lambda_t)$.

\hspace{.5cm}
\begin{myproof}{Theorem}{\ref{th:alpha,beta}}
Since $c\equiv 0\,(\text{mod }p)$, it follows from (\ref{w_p thm3}) that 
$$
W_p^{-1}\begin{pmatrix} a & b\lambda_t\\ c\lambda_t & d \end{pmatrix}W_p\in H_{M_p}(\lambda_t).
$$ 
Thus, $H_{M_p}$ is normalized by $W_p$ in $\operatorname{PSL}(2,\mathbb{R})$. The M\"obius transformation $W_p\tau=-\displaystyle\frac{1}{p\tau}$ induces an automorphism $\omega$ on $X_2=H_{M_p}\backslash\mathbb{H}$ such that the following diagram commutes:
\begin{center}
    \begin{tikzcd}
&\mathbb{H} \arrow[r, "W_p"] \arrow[d, "\pi_t'",labels=left]
& \mathbb{H} \arrow[d, "\pi_t'"] \\
&X_2 \arrow[r, "\omega",labels=below]
& X_2.
\end{tikzcd}
\end{center}
Moreover, by Lemma \ref{lem:Fricke}, $\omega:X_2\rightarrow X_2$ satisfies the following functional equations:
\begin{align*}
    &\varphi\circ\omega=1-\psi,\\
    &\psi\circ\omega=1-\varphi.
\end{align*}
Hence, for $z\in X_2$, we have $\varphi(\omega(z))=1-\psi(z)=1-\beta$ and $\psi(\omega(z))=1-\varphi(z)=1-\alpha$, i.e., $\omega$ interchanges $\alpha$ and $1-\beta$, and $\beta$ and $1-\alpha$. Thus, we deduce that $(1-\beta,1-\alpha)$ is also a solution to (\ref{mod eq}) and $P(1-\beta,1-\alpha)=0.$

\end{myproof}

\section*{Acknowledgements}

This article is a part of the author's doctoral research \cite{alam-thesis} under the guidance of Professor Toshiyuki Sugawa. The author would like to express his sincere thanks to Professor Toshiyuki Sugawa for proposing this topic and for valuable suggestions.

\end{document}